# Retour d'expérience interdisciplinaire : intégrer en Physique et Chimie des applications du cours d'Algèbre

Mickaël Bosco & Nicolas Michel

**Abstract.** As a Physicist and a Chemist, we have had the opportunity to teach Algebra in our engineering school in Aix-en-Provence. The article deals with the interactions made between Physics, Chemistry and Algebra.

**Keywords.** Algebra, Didactic, Innovation, Physics, Chemistry, Theory of Didactic Situations, Anthropological Theory of the Didactic, Fundamental Theorem of Algebra.

**Résumé.** Physicien et chimiste de formation, nous avons été amenés à enseigner de l'Algèbre en cycle préparatoire intégré au sein de l'école d'ingénieur ESAIP sur le campus d'Aix-en-Provence. Nos parcours respectifs nous ont permis d'illustrer les concepts et notions par des applications innovantes issues de nos domaines d'expertises respectifs. L'article portera donc sur un retour d'expériences croisées et son analyse réflexive.

**Mots-clés.** Algèbre, Didactique, Innovation, Physique, Chimie, Théorie des Situations Didactiques, Théorie Anthropologique du Didactique, Théorème fondamental de l'algèbre.

## Table des matières

## 1. Introduction

Les expressions algébriques représentent le quotidien des scientifiques. Quelle que soit leur discipline ou leur matière dominante, ils font un usage important du registre symbolique. Outre l'aspect purement calculatoire, ils utilisent ainsi de nombreuses règles sémantiques qui peuvent s'apparenter à une véritable linguistique.

L'étudiant, quant à lui, est graduellement exposé à ces langages. Depuis le lycée, à mesure que le degré de conceptualisation et d'abstraction augmente, les symboles mathématiques et lettres représentant des variables et des quantités, ou grandeurs physiques, sont de plus en plus présents jusqu'à devenir beaucoup plus fréquents que les chiffres.

Pour autant, les étudiants que nous rencontrons sont souvent démunis, confrontés à des difficultés lors de la manipulation de tels symboles et d'expressions les utilisant. En tant qu'enseignants, nous devons accompagner le développement de leur capacité à extraire une signification des symboles utilisés : représentent-ils un nombre, une quantité (auquel cas il faut, bien sûr, introduire les unités) ? Les opérateurs ont-ils un sens plus profond que l'opération qu'ils représentent (le symbole « = » dans une équation chimique représente ainsi la conservation de masse, la conservation des charges…) ? L'interprétation, bien sûr, dépend du contexte mais elle est sous-tendue par l'ancrage théorique et généraliste qu'apporte l'algèbre. En ce sens, et c'est l'objet de cet article, nous pensons qu'une approche pédagogique interdisciplinaire privilégiant les aller-retours entre théorie et application permet d'ancrer profondément ce sens du symbole au sens de Arcavi. La modernisation de l'enseignement passe par l'insertion multiple de situations concrètes de l'usage de l'algèbre. Par la suite, les étudiants vont pouvoir transposer cet usage à la résolution d'autres problèmes, et gagner ainsi en compétence.

Après avoir présenté le contexte de notre **retour d'expériences croisées** en Physique et Chimie, nous exposerons quelques situations représentatives de ce choix stratégique dans nos enseignements respectifs. Nous terminerons en présentant le bilan du ressenti des étudiants, et le nôtre, concernant l'acquisition des compétences visées.

## 2. Enseignement de l'Algèbre et besoin d'innovation

L'évolution des profils d'étudiants au cours des dernières décennies et des réformes successives ont amené les collèges, lycées, et également le supérieur à travailler sur une refonte des programmes et méthodes pédagogiques. Au sein de notre école d'ingénieur Esaip, une direction est consacrée à l'innovation en pédagogie et notamment aux moyens à mettre en place afin de diversifier les pratiques existantes.

Le niveau des classes est de plus en plus hétérogène, phénomène qui est également observable dans toutes les écoles du Supérieur d'après une étude menée par Le Prévost [9]. Nos étudiants peuvent accéder à l'école après le bac, au niveau de la classe préparatoire intégrée, en passant le concours Puissance Alpha ou par une admission parallèle s'ils ont préalablement suivi un autre cursus.

En fin de cursus, ils obtiennent un diplôme d'ingénieur, de niveau bac+5, dans l'une des deux spécialités suivantes : Ingénieur du Numérique et Ingénieur en Gestion des Risques et Environnement.

Les Mathématiques occupent une place importante au sein de nos formations puisqu'elles en représentent plus de 20 %. Différentes matières y sont enseignées telles que l'Analyse, l'Algèbre, la Géométrie, les Probabilités et le Traitement Statistique. Les contenus de ces matières sont régis par des syllabus présentant le programme semestriel de chacune, conçus de manière collective et revus chaque année, au cœur d'une démarche réflexive. Ces contenus précisent également les objectifs en termes de connaissances et compétences mais également les modalités d'évaluations.

À titre d'exemple, la matière Algèbre 1, dont le programme a été élaboré par le responsable des programmes scientifiques après échange avec les enseignants de ces matières au sein des différents campus, incluant les auteurs, qui a lieu au premier semestre de la première année, possède les caractéristiques suivantes :

| Algèbre 1 | | |
|---|---|---|
| Volume horaire encadré : **30,0 h** | Charge de travail : **54,0 h** | Coefficient matière : **3 / 10** |

### **Objectifs en matière de compétences :**

- S'engager dans une recherche, mettre en œuvre des stratégies : la découverte et l'exploitation de problématiques, la réflexion sur les démarches suivies, les hypothèses formulées et les méthodes de résolution

1. Manipuler les méthodes et les outils mathématiques permettant une modélisation : l'élaboration de liens déductifs ou inductifs entre différents éléments
2. Maîtriser les outils de base en algèbre : élaborer une représentation, changer de cadre, traduire des informations dans plusieurs registres.

## **Objectifs en matière de connaissances :**

- Étudier les structures algébriques usuelles et les liens existants entre elles ; propriétés de la divisibilité des entiers et des congruences.
- Étudier l'arithmétique des entiers relatifs et des polynômes à une indéterminée
- Assimilation des notions de base de l'algèbre linéaire
- Découvrir le caractère général des méthodes linéaires, notamment à travers leurs interventions en analyse et en géométrie.

## **Modalités pédagogiques :**

La pédagogie mise en œuvre développe la participation, la prise d'initiative et l'autonomie des étudiants. Plusieurs démarches peuvent être utilisées dans la scénarisation pédagogique avec l'introduction et le développement des outils numériques :

| Apprendre en mode de Cours/TD | *Présentations des notions nouvelles* |
|---|---|
| Apprendre en mode de cours inversé | *Consolidations des acquis* |

| Apprendre en mode échange | *Interagir « avec et entre » les élèves via Moodle :* <br> *forum, tests en ligne, etc.* <br> *Interagir à l'aide des boitiers de* |
|---|---|
| Apprendre en mode projet | *Mettre en pratiques des notions par approches* |
| Apprendre ensemble | *Organiser des séances de soutien et de tutorat entre les élèves* <br> *Elaborer des tests par les pairs* |

-
## Contenu pédagogique et progression :

1. <u>Structures algébriques</u> :
	- Concepts « lois de composition » interne, propriétés
	- Groupes, sous-groupes, morphismes
	- Anneaux, sous-anneaux, corps, anneaux intègres
	- Divisibilité, PGCD, lemme de Gauss, factorialité

2. <u>Polynômes et fonctions rationnelles</u> :
	- Concepts « Notations K[X] ».
	- Fonctions polynomiales, dérivées, racines, multiplicité
	- Division Euclidienne, algorithme d'Euclide pour le PGCD

3. <u>Systèmes & Matrices :</u>
	- Concepts et présentation des systèmes, types de résolution de base (substitution et combinaison), pivot de gauss et astuces.
	- Groupe ($M_{n,p}(K)$,+), anneau ($M_n(K)$,+,×), groupe ($Gl_n(K)$,×), systèmes de Cramer.
	- Cas de matrices triangulaires supérieures, Transposition et inversion.

- Calcul Matriciel : Lien avec les systèmes, notion de ligne et colonne, vocabulaire général, opérations sur les matrices, inversion de matrice.

-

| Séances | Chapitre / Activité |
|---|---|
| 1-6 | Structures Algébriques |
| 7-12 | Polynômes et fonctions rationnelles |
| 13-18 | Systèmes et matrices |

Veuillez noter que le contenu de la fiche matière de l'Algèbre du premier semestre est précisé ci-dessus, mais il est important de noter que la matière Algèbre est présente lors des trois premiers semestres de formation sur les 4 semestres du cycle préparatoire intégré. Comme cela est précisé précédemment, pour chaque matière, une fiche matière est créée selon le processus suivant. Dans un premier temps, les enseignants de la matière échangent sur le programme et une première proposition est réalisée et envoyée au responsable des programmes. Après avoir reçu l'ensemble des fiches matières, le responsable des programmes organise une réunion de travail avec l'ensemble des équipes enseignantes et des lectures croisées sont organisées afin de corriger et valider les contenus. A ce titre, nous avons donc participé à la construction des fiches matières d'Algèbre mais également de Physique et Chimie.

Les responsables des programmes et les équipes enseignantes de notre école travaillent actuellement sur la refonte en profondeur du programme pour les cinq années de formation. Cette mise à jour des contenus se révèle indispensable, imposée par les spécialités et options proposées au lycée depuis la mise en place de la réforme de 2019 : il nous faut prendre en compte la possibilité que des étudiants de première année n'aient pas tous suivis l'option « mathématiques expertes » dans le réajustement du programme en sciences. De fait, les nouveaux programmes se voient simplifiés au niveau des notions étudiées et recentrés en termes de contenus.

Ainsi, l'arrivée dans le supérieur des étudiants ayant connu les nouveaux programmes du lycée nous a amené à faire évoluer les fiches matières et notamment celles d'Algèbre selon les points suivants :

- Les structures algébriques (anneaux, corps, groupes, monoïdes, …) ne sont plus présentes en première année.

- La Décomposition en Eléments Simples, notée DES, n'est plus abordé dans cette matière.

- Le programme se focalise plus sur l'aspect calculatoire que sur l'aspect conceptuel.

- La finalité du cours est changée et recentrée sur les calculs matriciels.

Nous allons nous intéresser dans la suite à l'enseignement de la matière Algèbre et voir comment nos expériences en tant qu'enseignants en Physique et Chimie nous ont permis de mieux appréhender ces notions parfois abstraites pour les étudiants.

La diversité des parcours n'a fait qu'accentuer le besoin d'innover et ceci pour plusieurs raisons. La matière Algèbre n'est pas traitée au lycée en tant que matière à part entière. Or, ses concepts de base ne peuvent être réinventés, par conséquent, le programme doit maintenir l'apprentissage des chapitres clés en début de première année post-bac.

De plus, l'enseignement des notions théoriques, déjà compliqué en lui même, est rendu plus difficile encore par l'hétérogénéité du niveau des étudiants au sein d'une même classe.

Nos formations de Physicien et Chimiste n'étaient pas prédestinés à l'enseignement des Mathématiques à la base puisque nous sommes tous deux des enseignants issus d'un doctorat en Mécanique des Fluides et d'un doctorat en Chimie Organique donc plutôt orientés Physique et Chimie à la base. Même si ces cursus présentent une forte proportion de Mathématiques, elles étaient, le plus souvent, appliquées. Nous les avons donc appréhendées plus en tant qu'outils que finalité. Nos parcours professionnels respectifs nous ont cependant amenés à dispenser des enseignements en Mathématiques en parallèle de nos domaines de spécialité et ces parcours se sont révélés être une véritable force comme vous le verrez dans la partie suivante.

## 3. Applications de l'Algèbre en Physique

Nous avons enseigné, au sein de notre école, l'Algèbre et, dans le même temps, soit la Physique soit la Chimie au cours des deux années du cycle préparatoire. Ceci nous a permis de dresser un bilan puis, en prenant appui sur nos retours, de participer à l'élaboration des mises à jour des programmes.

Lorsque nous décortiquons les enseignements en Physique et Chimie, il se trouve que plusieurs notions abordées au sein de ces matières sont des applications directes de l'Algèbre.

Intéressons-nous tout d'abord à l'enseignement de Physique. Au sein de notre cycle préparatoire, la matière Physique, présente à chaque semestre, aborde différentes thématiques suivantes :

- Physique 1, semestre 1 : Optique géométrique (lois de Descartes, approximation de Gauss, application à la fibre optique, au prisme et à

l'arc-en-ciel) ; Mécanique du point (changement de repère, définitions, lois de Newton, Travail et énergie).

- Physique 2, semestre 2 : Mécanique du point et du solide (oscillateur harmonique, forces centrales, mécanique du solide, théorème du moment cinétique, gravitation).

- Physique 3, semestre 3 : Électromagnétisme (distributions discrètes et continues, théorème de Gauss, approche infinitésimale, symétries et invariances, potentiel électrostatique) ; Électrocinétique (associations de dipôles, modèle Thévenin-Norton, régime transitoire, circuits RC RL et RLC).

- Physique 4, semestre 4 : Thermodynamique physique (notions fondamentales, systèmes ouverts et fermés, premier principe de la Thermodynamique).

Les différents exemples cités pour les matière Physique et Chimie ont été utilisés dans nos cours de Physique 1, Physique 3 et Chimie 2.

Lors des enseignements effectués en Physique au cours des trois premiers semestres, l'utilisation des notions acquises en Algèbre est récurrente.

En effet, en fonction des notions abordées, la partie théorique algébrique viendra en amont et donc sera appliquée en Physique. On peut également découvrir une méthode en Physique et par la suite étudier la théorie algébrique derrière cette méthode. Ces notions sont intéressantes car : pour l'Algèbre, les étudiants découvrent une finalité des notions et, en même temps, pour la physique, ils découvrent des méthodes qui permettent de pallier leurs difficultés conceptuelles.

En effet, pour la matière Physique 1, les étudiants découvrent plusieurs repères d'étude : le repère cartésien, le repère polaire, le repère cylindrique et le repère sphérique. Une des premières difficultés des étudiants réside dans le passage d'un repère à un autre et cela peut se faire à l'aide de matrices de passage introduites en Algèbre comme nous l'avons illustré ci-dessous lors du passage du repère cartésien au repère polaire ou cylindrique [2] :

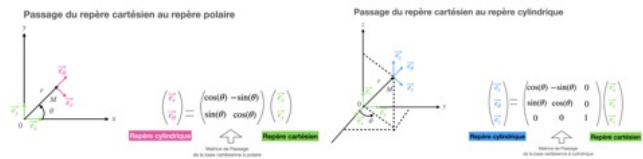

Figure 1 – Matrices de passages associées aux changements de repère

Ces matrices de passage permettent de changer de base et de se placer dans le repère souhaité dans un exercice de Physique d'une manière plus facile à utiliser que les projections tant redoutées des étudiants. De la même manière, le passage en repère sphérique nécessite une autre matrice de passage décrite dans la figure ci-après :

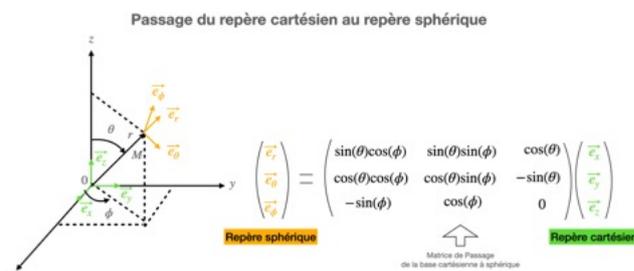

Figure 2 – Matrice de passage associée au changement de repère vers un repère sphérique

Ces matrices de changement de base sont introduites justement en cours d'Algèbre en première année.

Au cours de cet enseignement, ces matrices sont présentées comme des tableaux contenant des nombres comme les valeurs prises par une fonction trigonométrique, par exemple. L'approche qui consisterait à considérer les matrices à coefficients dans un ensemble de fonctions paraît trop compliquée en première année et n'est pas prévue dans le programme. Nous tenons à préciser que cette formule est démontrée aux étudiants, ce qui permet de leur expliquer l'origine des différents coefficients.

Les lois de Newton, introduites au lycée et revues en première année, sont également l'occasion d'utiliser des notions d'Algèbre comme notamment les équations différentielles et leurs systèmes d'équations associées.

La deuxième loi de Newton ci-dessous est une loi très importante en Mécanique puisqu'elle permet de décrire la dynamique d'un système en mouvement :

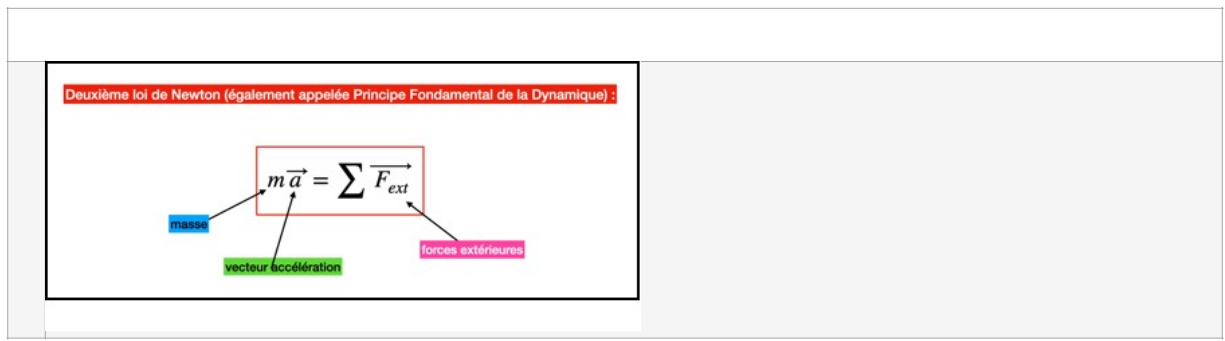

Figure 3 – Deuxième loi de Newton

En projetant cette relation vectorielle sur les différents axes, on peut être amené à résoudre un système d'équations différentielles comme celui ci-dessous qui correspond au modèle Lotka - Volterra [3] :

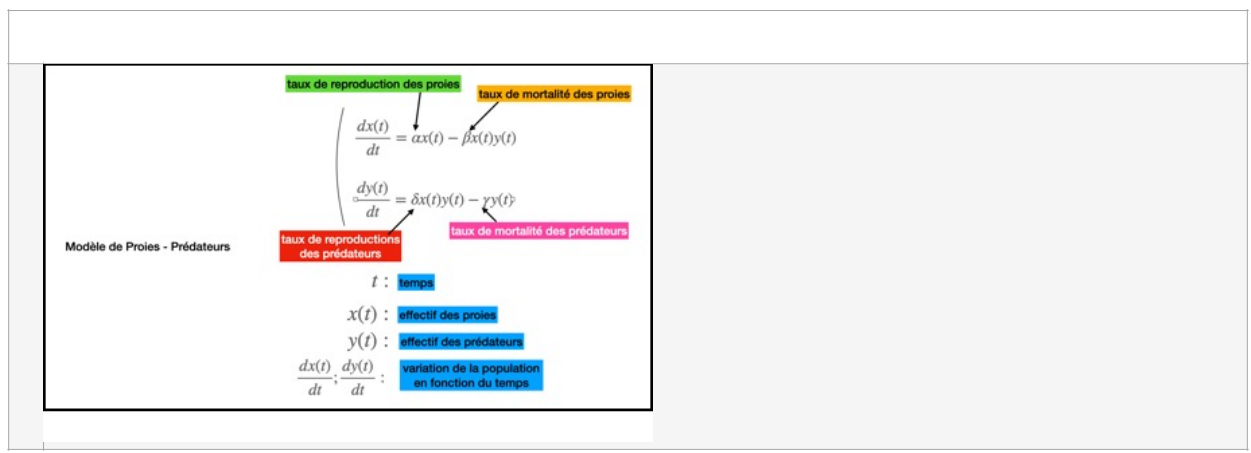

Figure 4 – Modèle proies-prédateurs

Ce système est étudié habituellement en Physique Non Linéaire pour aborder la théorie du chaos mais pour ma part j'ai l'habitude de le présenter aux étudiants dès le cours d'Algèbre 1. La résolution dépend fortement des conditions initiales. Il est préférable d'introduire cette étude de cas sous forme d'un jeu avec une mise en situation mêlant moutons et loups au sein d'un champ fermé, dont l'évolution temporelle va dépendre des conditions initiales.

Ces équations sont effectivement utilisées pour décrire la dynamique des systèmes biologiques dans lesquels un prédateur et sa proie interagissent, comme le lynx et le lièvre ou le loup et le mouton. Ce

système d'équations peut être résolu à l'aide d'équations différentielles comme cela est précisé ci-dessous :

$$\begin{cases} \dfrac{dx(t)}{dt} = \alpha x(t) - \beta x(t)y(t) \\ \dfrac{dy(t)}{dt} = \delta x(t)y(t) - \gamma y(t) \end{cases} \Rightarrow \begin{pmatrix} \dfrac{dx(t)}{dt} \\ \dfrac{dy(t)}{dt} \end{pmatrix} = \begin{pmatrix} \alpha & -\beta x(t) \\ \delta x(t) & -\gamma \end{pmatrix} \begin{pmatrix} x(t) \\ y(t) \end{pmatrix}$$

Vecteur dérivée — Matrice de Lotka Volterra — Vecteur position

Figure 5 – Système d'équations différentielles

En cours d'algèbre, le problème est présenté et discuté sous forme d'une étude de cas guidée dans le but d'introduire les systèmes d'équations différentielle. Par la suite, en Algèbre 3, on explique comment il est possible de résoudre ce système en réduisant la matrice de Lotka-Volterra pour étudier ses valeurs propres.

D'autres notions également étudiées en Physique 2 utilisent des notions d'Algèbre avec notamment l'oscillateur harmonique et les systèmes couplés, le théorème du moment cinétique et la mécanique du solide avec les matrices d'inertie.

Nous retrouvons ces dernières au sein des problèmes physiques dans lesquels nous ne pouvons pas réduire la dynamique d'un système à celle d'un seul point géométrique. Dans ce sens, une matrice d'inertie, dont les composantes sont écrites ci-dessus [4], représente la géométrie et la répartition de la masse du solide étudié.

Matrice d'inertie du solide S centrée au point O

$dm = dxdydz$ : Masse volumique

$$I(O,S) = \begin{pmatrix} \int_S (y^2 + z^2)dm & -\int_S yxdm & -\int_S zxdm \\ -\int_S xydm & \int_S (x^2 + z^2)dm & -\int_S zydm \\ -\int_S xzdm & -\int_S yzdm & \int_S (x^2 + y^2)dm \end{pmatrix}$$

$\int_S (y^2 + z^2)dm$ : Matrice d'inertie par rapport à (O,x)    $\int_S xydm$ : Matrice d'inertie par rapport à (O,x,y)

$\int_S (x^2 + z^2)dm$ : Matrice d'inertie par rapport à (O,y)    $\int_S yzdm$ : Matrice d'inertie par rapport à (O,y,z)

$\int_S (x^2 + y^2)dm$ : Matrice d'inertie par rapport à (O,z)    $\int_S xzdm$ : Matrice d'inertie par rapport à (O,x,z)

Figure 6 – Matrice d'inertie d'un solide

De plus, en Physique 3, pour les cours d'Électromagnétisme, différentes relations utilisent les dérivées vectorielles comme pour le potentiel vecteur par exemple. D'autres opérateurs différentiels tels que la Divergence, le Rotationnel et le Laplacien nécessitent également les dérivées vectorielles et donc une approche matricielle est privilégiée [5].

L'ensemble des exemples ci-dessus sont abordés dans le cours de Physique depuis plusieurs années par l'enseignant de Physique et par la suite les notions algébriques sont également revues en cours d'Algèbre en tant qu'approfondissement.

En Électricité, les études des dipôles RC, RL et RLC utilisent des équations différentielles dont certaines peuvent être couplées et résolues à l'aide de matrices comme nous avons pu le voir précédemment avec le système différentiel. On remarque ainsi que plusieurs notions algébriques sont utilisées en Physique et nous allons voir par la suite que c'est également le cas en Chimie. Dans la suite, nous nous intéresserons aux applications de l'Algèbre en Chimie.

## 4. Applications de l'Algèbre en Chimie

En Chimie, le programme est orienté vers l'utilisation des modèles. Dès le début de la formation, les différents modèles atomiques sont présentés en insistant sur leur évolution historique et théorique, leurs avantages et inconvénients. La progression employée a pour but, outre l'apport au niveau de la culture scientifique, de développer l'esprit critique des étudiants vis à vis des outils à leur disposition, de leur montrer que les sciences expérimentales sont en perpétuel développement et que les avancées théoriques et technologiques permettent d'affiner les hypothèses et d'élaborer des modèles de plus en plus précis ou adaptés. S'en suivent les transformations physiques et chimiques de la matière, jusqu'au bilan de matière et aux équilibres avec une initiation à la thermochimie et à la cinétique. Par la nature même des phénomènes étudiés et leurs échelles, les étudiants se rendent très rapidement compte de l'utilité des modèles. Sans aller jusqu'au cycle inductif-déductif permettant leur élaboration, nous les accompagnons dans la phase d'identification du modèle à utiliser dans un contexte donné. Cette approche nous permet de contextualiser le modèle en tant qu'outil et de faire sens à travers la modélisation.

Or, nous parlons ici de modélisation algébrique des phénomènes : un modèle est, par essence, par construction, algébrique ! Les points de contact entre disciplines sont alors multiples et les applications pratiques de l'algèbre en chimie foisonnent avec de nombreuses références théoriques mais aussi sur la didactique associée : stœchiométrie, cinétique chimique, chimie combinatoire, chimie des polymères, topologie et descripteurs moléculaires, isomérie, stéréochimie, problèmes d'optimisation, … pour n'en citer que quelques-unes.

Nous nous contenterons, dans cet article, de quelques exemples représentatifs tirés des cours de Chimie et d'Algèbre dispensé au premier et second semestre de la première année du cycle préparatoire intégré. L'exemple le plus marquant ici est la résolution matricielle appliquée à l'ajustement des coefficients stœchiométriques d'une équation chimique. Il s'agit également de l'exemple qui a reçu la plus forte adhésion cette année. Il se trouve en effet que l'étude des transformations chimiques coïncide avec l'étude des matrices en tant qu'outils de résolution de systèmes d'équations linéaires en algèbre linéaire.

Outre la résolution de systèmes, cet exemple permet également de revenir sur les conventions d'écriture symbolique des espèces chimiques, ions, ions polyatomiques et molécules, qui constituent un langage à part entière [6]. En effet, la ligne d'information que représente l'écriture d'une transformation chimique contient une multitude d'informations codées dont le sens sous-jacent est loin d'être évident pour nos étudiants avec des biais cognitifs et des difficultés de compréhension/appropriation :

• Difficulté de préemption des échelles puisqu'une formule représente à la fois une entité chimique et un grand nombre de ces entités (microscopique – entité chimique, macroscopique – quantité de matière)

• Difficulté d'interprétation des réorganisations chimiques, ruptures ou recréation de liaisons, et physiques à savoir ma modification d'interactions amenant à des états physiques différents.

Ainsi, exemple d'utilisation de l'algèbre dans le cours de chimie, la symbolique de l'entité chimique polyatomique peut s'appréhender comme représentation d'une combinaison linéaire, en considérant les atomes ou les grandeurs associées à ces atomes comme vecteurs et leur nombre comme des scalaires. Le premier exemple donné par l'un des auteurs, durant les séances de travaux encadrés consiste en un calcul de masse molaire moléculaire, somme des masses molaires atomiques des éléments chimiques constituant l'entité chimique, présenté explicitement comme une combinaison linéaire : On donne la formule brute du glucose, calculer sa masse molaire moléculaire en utilisant les données issues de la classification périodique à votre disposition. Il est demandé oralement aux étudiants de retrouver l'équation suivante en insistant sur l'identification de ses constituants :

Formule brute du glucose :

Équation :

En obligeant les étudiants à formaliser leur raisonnement, cette approche permet d'asseoir le concept même de formule brute (sa signification et son utilisation). Elle prend tout son sens, après cette première question, lorsque le problème est de la retrouver à partir de grandeurs liées aux atomes. Par exemple :

Connaissant la masse molaire moléculaire et la composition centésimale massique du glucose, 40,00 % de carbone, 6,67 % d'hydrogène et 53,33 % d'oxygène, retrouver sa formule brute.

La consigne et son traitement sont alors clairs et évitent toute ambiguïté, interprétation ou raisonnement alternatifs. La grande majorité des étudiants confrontés à cet exercice produit le raisonnement suivant sans hésitation :

Pour

La seconde partie présente la combustion complète du glucose et l'exercice porte alors sur l'ajustement des coefficients stœchiométriques de l'équation de cette transformation : « Du glucose réagit avec du dioxygène pour donner de l'eau et du dioxyde de carbone. Écrire l'équation associée à cette réaction. ». On cherche, ici, dans quelles proportions les entités chimiques doivent réagir, pour les réactifs, et apparaître, pour les produits, pour que les principes de conservation soient vérifiés.

Au lycée, dès la seconde, avec des notions telles que la conservation de masse abordées dès le collège, l'approche d'un tel exercice est majoritairement déductive et intuitive, et, vraiment minoritairement, algorithmique et algébrique. Interrogés sur leur stratégie, la réponse la plus fréquente des étudiants reste : « on compte le nombre d'atomes de chaque élément à droite et à gauche de l'équation et on ajuste les coefficients si besoin ».

Bien qu'exposée clairement, cette stratégie n'est pas claire pour eux et ils peinent encore à ajuster les coefficients même dans des cas simples, comme celui proposé, en raison des difficultés et biais reconnus suivants :

- Difficulté à transposer le nom des entités en symboles chimiques
- Difficultés au niveau de la composition de ces entités

- Difficulté dans le choix de l'élément à considérer « en premier »
- Difficulté à traiter correctement les charges

Le parti pris a été de présenter la méthode matricielle dès le départ. Le cours d'Algèbre, mené par le même auteur, sur l'élimination de Gauss ayant eu lieu peu avant, la chronologie a été parfaite.

Nous avons donc transformé l'équation de réaction, en introduisant les coefficients/inconnues, en un système d'équations linéaires à partir de chacun des éléments présents dans les entités moléculaires et construit la matrice correspondante. Ce processus a ainsi permis de tester la compréhension de l'étape précédente et l'appropriation par les étudiants du concept de combinaison linéaire pour une entité chimique. Du point de vue de l'Algèbre, il a également permis de revenir sur le nombre de solutions d'un système d'équations linéaire, la construction de la matrice associée à un tel système et la relation AX=B, avec A matrice des coefficients, X matrice des variables et B matrice des solutions.

La méthode présentée aux étudiants est la suivante : les coefficients stœchiométriques cherchés sont désignés par des lettres, représentant autant d'inconnues (ce qui invite fortement les étudiants à une résolution algébrique). On écrit alors l'équation bilan de réaction :

Le système se construit rapidement en faisant un bilan de chaque élément :
- Pour le carbone, C, le réactif glucose comporte 6 atomes de carbone mais ce nombre se trouve multiplié, dans l'équation de réaction, par le coefficient stœchiométrique , soit . Au niveau des produits, le carbone ne se retrouve que dans le dioxyde de carbone (1 atome), multiplié là par le coefficient . Ainsi,

- Pour l'hydrogène, H,

• Pour l'oxygène, O,

Soit, en insistant sur le fait que ces équations doivent être résolues simultanément, le système suivant par identification :

Ce système comporte trois lignes d'information, trois équations, pour quatre inconnues et constitue ainsi un système sous-déterminé avec une infinité de solutions. Ceci a permis l'évocation du rang, encore non étudié au moment de l'exercice, en constatant l'indépendance linéaire des trois équations.

Remarque : les étudiants s'attendant à une solution unique, ce résultat leur a paru contre-intuitif et un rappel de la définition des coefficients stœchiométriques en tant que coefficients de proportionnalité a été nécessaire pour l'expliquer.

En ce qui concerne la résolution, nous avions plusieurs possibilités mais la construction d'une matrice carrée 3x3 et d'une matrice des coefficients dépendants de la dernière inconnue a été retenue (sous la forme $AX=B\delta$) car la plus proche de ce qui avait été traité en algèbre :

Remarque : la matrice A est, ainsi, la matrice des vecteurs de $\mathbb{R}^3$ correspondant au nombre de chaque atome dans chacune des entités chimiques, réactifs (positifs) et produits (négatifs) hormis le dernier.

Soit

pour les réactifs et

pour le produits. D'où

La matrice B correspond, elle, au dioxyde de carbone :

Et la matrice augmentée, A augmentée de B, sur laquelle les étudiants ont appliqué l'algorithme d'élimination de Gauss présenté durant le cours d'algèbre :

En retournant au système simplifié, ils se sont ainsi rendu compte que les solutions dépendaient de la valeur de $\delta$, ce qui a conforté les explications sur l'infinité de solutions.

Il est généralement préférable d'obtenir des coefficients entiers donc, en posant $\delta = 6$, les inconnues $\alpha$, $\beta$ et $\gamma$ peuvent être calculées et l'équation de réaction devient :

Après quelques exercices sur le même modèle, accompagnés, en autonomie, avec correction par l'enseignant puis par les pairs, les équations contenant des ions, donc des charges ont été introduites. Il s'agissait là d'utiliser la conservation de la charge (la somme des charges doit être égale de part et d'autre de l'équation - notion admise par l'ensemble des étudiants), sous forme d'une combinaison linéaire à déterminer. Le premier exemple traité a été la réaction du thiosulfate sur le Diiode donnant les ions iodure et Tétrathionate :

Laissés en autonomie pour rechercher une solution, les étudiants ont suivi le même raisonnement que précédemment avec les atomes, ce qui les a rapidement conduits à l'équation associée aux charges :

La résolution complète n'a alors posé aucun problème.

Ainsi les étudiants se sont approprié assez rapidement la méthode et ont pu l'appliquer à des cas plus complexes (avec l'introduction des charges mais aussi avec des réactions moins évidentes comme

 [7]).

Ils ont également vu la possibilité d'automatisation de la tâche en utilisant l'algorithme associé à la méthode d'élimination de Gauss. Cette méthodologie pour équilibrer les équations a été traité en classe depuis plusieurs années également.

## 5. Bilan et conclusion

Les résultats de ce retour d'expériences croisées sont présentés ci-après, sous forme d'un bilan. En effet, les différents exemples, méthodes et applications, présentés au sein de nos enseignements, ont donné lieu à un questionnement des étudiants dans une approche réflexive : depuis que nous enseignons au sein de ce cycle, nous avons pu échanger avec les étudiants en chaque fin de semestre sur leur ressenti à ce sujet.

Afin d'illustrer les retours étudiants, nous réalisons à chaque fin de semestre des enquêtes qualité pour faire le bilan sur chaque matière enseignée et mettre en avant les corrections à apporter pour l'année suivante ou au contre péréniser les pratiques qui ont eu un bon retour. Ce questionnaire est réalisée en classe entière lors d'un créneau fixé dans l'emploi du temps afin d'avoir un taux de réponse de 100%.

Dans ce sens, vous trouverez ci-dessous un extrait des questions concernant notre étude précédente qui est réalisée pour chaque matière et chaque semestre, ici pour la matière Algèbre 1 à titre d'exemple :

Question 1 – Quelle est votre appréciation sur le contenu des cours ( supports, interactivité, cohérence cours, TD/TP/Projets, ..) de la matière Algèbre 1 ?

Question 2 – Quelle est votre appréciation sur les évaluations réalisées ( types, modalités, ) de la matière Algèbre 1 ?

Question 3 – Quelle est votre appréciation sur les méthodes pédagogiques utilisées de la matière Algèbre 1 ?

Question 4 – Avez-vous des commentaires libres à formuler pour la matière Algèbre 1 ?

Pour notre étude, cette enquête a été renseignée par une promotion de 60 étudiants sur les 4 semestres du cycle préparatoire intégrée et les résultats obtenus ont mis en avant les points suivants :

En analysant les réponses aux trois premières questions couplées aux commentaires libres, nous avons pu mettre en avant le fait que dans la grande majorité, les étudiants sont favorables à ces approches et ces sondages nous informent que plus de 80% des étudiants en sont ravis.

Ils se rendent compte que les compétences liées aux Mathématiques et Sciences Expérimentales sont liées. Les étudiants de première année, issus d'un bac général ou technologique, découvrant l'Algèbre, sont souvent déconcertés par le côté abstrait et nouveau de ce langage qui est loin d'être intuitif pour eux [8]. Notre objectif est donc de leur fournir des applications, utilisations… permettant de les familiariser avec les symboles, notions et concepts dans un cadre différent du cadre formel, sans théorème ni démonstration, pour redonner du sens. Lors de l'utilisation de ces notions dans des exercices ciblés de physique ou de chimie, ils s'aperçoivent qu'ils font de l'Algèbre depuis longtemps, manipulation de formules, et qu'en définitive, dans l'ensemble, ils maîtrisent plutôt bien ces notions.

Ceci a également un impact sur leurs résultats puisque le taux de réussite, lors d'évaluations sur les notions traitées de manière interdisciplinaire, augmente à hauteur moyenne de 30 % ce qui n'est pas négligeable en sciences, évolution de la moyenne de la classe sur une évaluation type, similaire en contenus, avant et après utilisation interdisciplinaire, en physique et en chimie. Cette hausse est observable dans la validation des différentes unités d'enseignements au cours des deux premières années de formation et plus précisément des matières Physique, Chimie et Algèbre. De plus, il y a eu un effet notoire en cours d'Algèbre une fois que ces notions ont été utilisées en sciences physiques car les étudiants se rappellent les éléments mathématiques qui ont servi à l'élaboration des cours en Physique et Chimie. Cette augmentation du taux de réussite est donc mesurable au sein des matières concernées, à savoir Algèbre, Physique et Chimie sur chaque semestre et donc par conséquent a entrainé une augmentation du taux de réussite à l'échelle de chaque semestre et de l'année académique.

Outre l'augmentation du taux de réussite des apprenants qui est quantifiable, nous avons pu nous rendre compte que l'aisance des étudiants concernant certaines notions théoriques, comme par exemple ici les notions de changement de base par le biais des matrices des passage est observable. En effet, même si dans un premier temps, en cours d'Algèbre

1, les étudiants se contentent d'apprendre par cœur cette notion, le fait de travailler de nouveau ce thème dans une autre matière, comme ici en Physique 2 ou en Chimie 3, va permettre à l'apprenant de mieux comprendre le sens théorique qui correspond dans le cas des matrices de passage au sens géométrique et donc de mieux retenir sur le long terme.

Pour l'enseignant, outre les points soulevés précédemment, l'approche interdisciplinaire amène les avantages suivants : elle permet tout d'abord de diversifier les méthodes d'enseignement en multipliant les cas d'utilisation hors du cadre du cours. Ainsi, l'enseignant peut choisir des problèmes ou exemples pertinents, se prêtant bien à l'observation des notions et leur intégration dans un contexte motivant et faisant sens pour les étudiants. La construction des séances pédagogiques en est d'autant facilitée avec la possibilité qui en découle de mettre en place un apprentissage par résolution de problèmes. Elle permet également de lever certains biais de compréhension.

De plus, la mise en place de cette interdisciplinarité peut permettre à l'enseignant de renforcer les liens existants entre les différentes matières scientifiques et peut lui permettre de gagner du temps dans sa séquence pédagogique. Lorsqu'une notion théorique a été préalablement étudiée dans une autre matière, si l'apprenant étudie de nouveau cette notion, dans le même contexte ou non, celui-ci va pouvoir la comprendre plus rapidement, ce qui peut être bien utile lorsque le volume horaire d'une matière est restreinte par rapport aux compétences et connaissances à acquérir.

En revanche, et pour qu'elle soit pleinement efficace, l'approche doit impérativement être présentée et explicitée en amont aux étudiants. Le but n'est pas d'exposer l'Algèbre comme une discipline « omnipotente » ou incontournable. Ceux qui présentent des réticences ou qui ne se sentent « pas bons » dans cette discipline ne doivent pas se sentir exclus. En cela, et comme précisé plus haut, les exemples et problèmes se doivent d'être judicieusement choisis.

Il est important également de souligner qu'étant donné que le cours d'Algèbre est présent lors des trois premiers semestres de formation et que les cours de Physique et Chimie sont présents lors de l'ensemble des quatres semestres du cycle préparatoire, la mise en place de cette interdisciplinarité est bénéfique à la fois pour les sciences physiques mais également pour l'Algèbre. En effet, les étudiants découvrent l'utilité de cette matière abstraite dans un autre domaine et, par la suite, une fois qu'ils l'ont utilisé dans une autre thématique, ils mémorisent mieux et comprennent mieux les techniques. Par conséquent cette interdisciplinarité favorise la cohésion de ces matières mais également leurs apprentissages croisées.

Avec un recul de maintenant plus de quatre ans, les bilans du côté étudiant comme enseignant se révèlent largement positifs : pour les étudiants, l'interdisciplinarité permet une montée en compétences, alors que pour les enseignants cela permet de diversifier les méthodes didactiques et d'innover.

De plus, avant d'avoir eu la possibilité d'enseigner ces deux matières, même en tant qu'enseignant, nous ne nous rendions pas spécialement compte de l'interdisciplinarité et nous menions nos enseignements de manière indépendantes sans montrer les différents liens qui existent.

Auparavant, avant d'enseigner au sein de notre école d'ingénieur actuelle, nous n'intervenions pas sur plusieurs matières d'enseigner différentes, nous enseignions soit uniquement les mathématiques soit les sciences physiques et par conséquent nous n'avions pas le recul nécessaire pour pratique cette interdisciplinarité. Ce que nous avons pu observer cependant, mais pas quantifier, c'est que les étudiants acquiert plus rapidement les connaissances et compétences par le biais de cette interdisciplinarité.

Cet article montre que les approches transverses avec l'Algèbre comme discipline commune permettent une montée en compétences des étudiants : ils peuvent réinvestir leurs connaissances dans d'autres contextes que celui du cours de Mathématiques. Nous sommes persuadés que leur montrer cette transversalité, et ainsi leur prouver l'universalité du langage et des outils algébriques, peut les conduire à utiliser ces outils et dans des environnements totalement différents. L'objectif étant, bien évidemment, de les amener à manipuler et produire eux-mêmes des expressions algébriques modélisant les problèmes qu'ils rencontrent, de les interpréter et d'utiliser les outils adaptés à leur résolution. Le champ d'application peut alors s'élargir : informatique, statistiques, économie, etc.

Comme présenté ci-dessus, l'Algèbre est une discipline fondamentale dispensée, en parallèle d'autres matières scientifiques telles que la Physique et la Chimie. Enseignant ces différentes matières, nous avons pu proposer des méthodes, exemples et applications interdisciplinaires qui ont conduit à une montée en compétences de nos étudiants. Cette approche pluridisciplinaire représente à elle seule une évolution des pratiques didactiques. Nous sommes bien conscients du fait que la situation présentée ici est particulière et que les enseignants se cantonnent généralement à leur domaine mais les avantages indéniables de cette approche nous conduisent à proposer une réflexion sur l'élaboration parallèle des déroulés pédagogiques et la concertation, indispensable, entre les enseignants des différentes disciplines.

L'interdisciplinarité est de plus en plus omniprésente dans l'enseignement et notamment en sciences. Cette pratique est particulièrement intéressante pour les notions qui peuvent paraître abstraite car elle permet de rendre plus concret des notions théoriques et ainsi être assimilées plus facilement pour les étudiants. De plus, les outils mathématiques font parti d'un bagage théorique nécessaire et obligatoire que doit avoir tout étudiant diplômé d'une école d'ingénieur de part l'omniprésence des mathématiques dans son futur métier. Il est donc important de montrer aux étudiants en quoi les mathétmatiques peuvent également être présents dans d'autres domaines et dans les autres enseignements de la formation suivie.

# Références

Mickaël Bosco
Enseignant chercheur à l'ESAIP, membre du groupe PION, IRES Marseille, CERADE
*e-mail:* mbosco@esaip.org

Nicolas Michel
Enseignant chercheur à l'ESAIP, CERADE
*e-mail:* nmichel@esaip.org